\documentclass[12pt,leqno]{article}
\usepackage{amsfonts}
\usepackage{amsmath, amsthm, amsfonts, amssymb, color}
\usepackage{mathrsfs}
\usepackage{color}
\theoremstyle{definition}

\newcommand{\scr}[1]{\mathscr #1}
\definecolor{wco}{rgb}{0.5,0.2,0.3}

\numberwithin{equation}{section} \theoremstyle{remark}

\newcommand{\ua}{\uparrow}

\title{{\bf  Integration by Parts Formula and Applications  for SPDEs with Jumps}\footnote{Supported in
 part by  NNSFC (11131003, 11431014), the 985 project, the Laboratory of Mathematical and  Complex Systems.} }
\author{
{\bf     Feng-Yu Wang  }\\
\footnotesize{ School of Mathematical Sciences,
Beijing Normal
University, Beijing 100875, China}\\
 \footnotesize{ Department of Mathematics,
Swansea University, Singleton Park, SA2 8PP, United Kingdom}\\
\footnotesize{  wangfy@bnu.edu.cn, F.-Y.Wang@swansea.ac.uk}}
\begin{document}
\allowdisplaybreaks
\def\R{\mathbb R}  \def\ff{\frac} \def\ss{\sqrt} \def\B{\mathbf
B}
\def\N{\mathbb N} \def\kk{\kappa} \def\m{{\bf m}}
\def\ee{\varepsilon}\def\ddd{D^*}
\def\dd{\delta} \def\DD{\Delta} \def\vv{\varepsilon} \def\rr{\rho}
\def\<{\langle} \def\>{\rangle} \def\GG{\Gamma} \def\gg{\gamma}
  \def\nn{\nabla} \def\pp{\partial} \def\E{\mathbb E}
\def\d{\text{\rm{d}}} \def\bb{\beta} \def\aa{\alpha} \def\D{\scr D}
  \def\si{\sigma} \def\ess{\text{\rm{ess}}}
\def\beg{\begin} \def\beq{\begin{equation}}  \def\F{\scr F}
\def\Ric{\text{\rm{Ric}}} \def\Hess{\text{\rm{Hess}}}
\def\e{\text{\rm{e}}} \def\ua{\underline a} \def\OO{\Omega}  \def\oo{\omega}
 \def\tt{\tilde} \def\Ric{\text{\rm{Ric}}}
\def\cut{\text{\rm{cut}}} \def\P{\mathbb P} \def\ifn{I_n(f^{\bigotimes n})}
\def\C{\scr C}      \def\aaa{\mathbf{r}}     \def\r{r}
\def\gap{\text{\rm{gap}}} \def\prr{\pi_{{\bf m},\varrho}}  \def\r{\mathbf r}
\def\Z{\mathbb Z} \def\vrr{\varrho} \def\ll{\lambda}
\def\L{\scr L}\def\Tt{\tt} \def\TT{\tt}\def\II{\mathbb I}
\def\i{{\rm in}}\def\Sect{{\rm Sect}}  \def\H{\mathbb H}
\def\M{\scr M}\def\Q{\mathbb Q} \def\texto{\text{o}} \def\LL{\Lambda}
\def\Rank{{\rm Rank}} \def\B{\scr B} \def\i{{\rm i}} \def\HR{\hat{\R}^d}
\def\to{\rightarrow}\def\l{\ell}\def\iint{\int}
\def\EE{\scr E}
\def\A{\scr A}
\def\BB{\scr B}\def\Ent{{\rm Ent}}

\maketitle

\begin{abstract}  By using the Malliavin calculus and finite jump approximations, the Driver-type integration by parts formula is established for the semigroup associated to stochastic (partial) differential equations with  noises containing a subordinate Brownian motion. As applications, the shift Harnack inequality and heat kernel estimates are derived. The main results are illustrated by SDEs driven by $\aa$-stable like processes.
\end{abstract} \noindent
 \vskip 2cm

\section{Introduction}

A significant application of the Malliavin calculus is to describe the density of a Wiener functional using the integration by parts formula.
In 1997, Driver \cite{Driver} established the following integration by parts formula for the heat semigroup $P_t$ on a compact Riemannian manifold $M$:
$$P_t (\nn_Z f)= \E\{f(X_t)N_t\},\ \ f\in C^1(M), Z\in \scr X,$$ where $\scr X$ is the set of all smooth vector fields on $M$, and $N_t$ is a random variable depending on $Z$ and the curvature tensor. From this formula we are able to characterize   the derivative w.r.t. the second variable $y$ of the heat kernel $p_t(x,y)$, see \cite{W12} for a recent study on integration by parts formulas and applications for SDEs/ SPDEs   driven by Wiener processes. The backward coupling method developed in \cite{W12} has been also used in \cite{Fan,Zhang} for SDEs driven by fractional Brownian motions and SPDEs driven by Wiener processes.  The purpose of this paper is to investigate the integration by parts formula and applications for SDEs/SPDEs  driven by purely jump L\'evy noises, in particular, to derive estimates on the heat kernel and its derivatives for the solutions.

Let $\scr L(\H)$ denotes the class of all bounded linear operators on $\H$ equipped with the operator norm $\|\cdot\|$. Let  $$\si: [0,\infty)\to\scr L(\H),\ \ b: [0,\infty)\times \H\to \H$$ be measurable and  locally bounded, such that $b_t:\H\to \H$ is Lipschtiz continuous locally uniformly in $t$. Consider the following stochastic equation on a separable Hilbert space $\H$:
\beq\label{1.1} X_t= \e^{At} X_0+ \int_0^t \e^{A(t-s)} b_s(X_s)\d s+\int_0^t \e^{A(t-s)} \si_s \d W_{S(s)} +V_t,\ \ t\ge 0,\end{equation}
where
 $W:=(W_t)_{t\ge 0}, S:=(S(t))_{t\ge 0}$ and $V:=(V_t)_{t\ge 0}$ are independent stochastic processes such that
\beg{enumerate} \item [(i)] $W$ is the cylindrical Brownian motion on $\H$ with $W_0=0$;
\item[(ii)]  $V$ is a c\'adl\'ag
  process on $\H$ with $V_0=0$;
  \item[(iii)] $S$ is the subordinator induced by a Bernstein function $B$, i.e. $S$ is a  one-dimensional   increasing L\'evy process with $S(0)=0$ and Laplace transform
$$\E \e^{-r S(t)}= \e^{-t B(r)},\,\,\ \ t,r\ge 0.$$
\item[(iv)] $(A,\D(A))$ is a linear operator   generating a $C_0$ contraction semigroup $\e^{At}$ on $\H$ such that
$$\int_0^T \|\e^{At}\|_{HS}^2\d t<\infty,\ \ T>0,$$ where $\|\cdot\|_{HS}$ is the Hilbert-Schmidt norm.
 \end{enumerate}
Then $(W_{S(t)})_{t\ge 0}$ is a  L\'evy process known as the subordinate Brownian motion (or subordinated process of the Brownian motion) with subordinator $S$ (see e.g.  \cite{A,J}),  and for any initial value $X_0=x\in \H$ the equation \eqref{1.1} has a unique solution (see \cite[Proposition 4.1]{WW}). Let $P_t$ be the associated Markov operator, i.e.
$$P_t f(x)= \E f(X_t(x)),\ \ f\in\B_b(\H), t\ge 0, x\in\H.$$
Bismut formula and Harnack inequalities for $P_t$  have been studied in \cite{Z} and \cite{WW} by using regularization approximations of $S(t)$, but the study of the integration by parts formula and shift Harnack inequality is not yet done.

\

 Since Ker$(\e^{At})=\{0\}$, the inverse operator $\e^{-At}: {\rm Im} (\e^{At})\to\H$ is well defined. To establish the integration by parts formula, we need the following assumptions.

\beg{enumerate} \item[(H1)]   $b_t\in C^2(\H)$, and for $(t,x)\in [0,\infty)\times \H$  there holds $\nn b_t(x): {\rm Im} (\e^{At})\to {\rm Im} (\e^{At})$   such that
$B_t(\cdot):= \e^{-At} (\nn b_t(\cdot)) \e^{At}$ satisfies
$$\|B_t\|_\infty:=\sup_{x\in\H} \|B_t(x)\|\le K_1(t),\ \ \|\nn  B_t(x)\|_\infty:=\sup_{x\in\H}\| \nn  B_t(x)\|\le K_2(t),\ \ t\ge 0$$ for some increasing $K_1,K_2\in C([0,\infty))$.
\item[(H2)] $\si_t$ is invertible   such that for some increasing $\ll_1,\ll_2\in C([0,\infty)),$
$$\|\si_t\|\le \ll_1(t),\ \ \|\si_t^{-1}\|\le \ll_2(t),\ \ t\ge 0.$$
\end{enumerate}

(H2) is a standard non-degenerate assumption, while (H1) means that the interaction between far away directions are weak enough. For instance, 
letting $-A$ be self-adjoint with discrete eigenvalues $0< \ll_1\le \ll_2\le\ll_3\cdots$ and eigenbasis $\{e_i\}_{i\ge 1}$, (H1) holds provided 
$$|\<\nn_{e_i} b_t, e_j\>|\le K_1(t)\e^{-t|\ll_i-\ll_j|},\ \ \|\nn\<\nn_{e_i} b_t, e_j\>\| \le K_2(t) \e^{-t|\ll_i-\ll_j|},\ \ t\ge 0, i,j\ge 1.$$

As already observed in \cite{W12} that comparing with the Bismut formula,   the integration by parts formula is usually harder to establish.
 To strengthen this observation, we explain below   that the regularization argument used in \cite{Z} for the Bismut formula is no longer valid for the integration by parts formula.
For simplicity, let us consider the case when $\H=\R^d, A=0, V_t=0, b_t=b$ and $\si_t=\si$.  As in \cite{Z}, for any $\vv>0$ let
$$S_\vv(t)= \ff 1 \vv\int_t^{t+\vv} S(s)\d s+\vv t,\ \ t\ge 0.$$ Then $S_\vv(\cdot)$ is differentiable and $S_\vv\downarrow S$ as $\vv\downarrow 0.$ Consider the equation (note that we have assumed $V_t=0, b_t=b$ and $\si_t=\si$)
$$\d X_t^\vv= b(X_t^\vv)\d t +\si \d W_{S_\vv(t)},\ \ X_0^\vv = X_0.$$
To apply the existing derivative formulas for SDEs driven by the Brownian motion, we take
$Y_t^\vv= X_{S_\vv^{-1}(t)}^\vv$ so that this equation reduces to
$$\d Y_t^\vv= b (Y_t^\vv) (S_\vv^{-1})'(t) \d t +\si \d W_t,\ \ Y_0^\vv=X_0.$$ In \cite{Z}, by using a known Bismut formula for $Y_t^\vv$ and letting $\vv\to 0$, the corresponding formula for $X_t$ is established. The crucial point for this argument is that the Bismut formula for $Y_t^\vv$   converges as $\vv\to 0$. However, since $S$ is not differentiable,   the existing integration by parts formula of $Y_t^\vv$ (see e.g. \cite[Theorem 5.1]{W12} with $H=\R^d$ and $ A=0$)
 $$\E(\nn_v f)(Y_T^\vv) =\ff 1 T \E\bigg\{f(Y_T^\vv)\int_0^T\big\<\si^{-1}\{v-t(S_\vv^{-1})'(t)\nn_v b(Y_t^\vv)\}, \d W_t\big\>\bigg\}$$
   does not converge to any explicit formula as $\vv\to 0$, except when  $\nn_v b$ is trivial.

So, to establish the integration by parts formula, we will take a different approximation argument, i.e. the finite jump approximation   used in \cite{WXZ} to establish the Bismut formula for   SDEs with  multiplicative  L\'evy noises.    We have to indicate that in this paper we are not able to establish the integration by parts formula for SDEs with multiplicative  L\'evy noises. Note that even for SDEs driven by multiplicative Gaussian noises, the existing integration by parts formula using the Malliavin covariant matrix is in general less explicit.

\

 To state our main result, for any $s\ge 0$ we introduce the $\scr L(\H)$-valued processes $(J_{s,t})_{t\ge s}$ and $(\tt J_{s,t})_{t\ge s}$, which solve  the following random ODEs:
 \beq\label{*0}  \ff{\d}{\d t} J_{s,t}= B_t(X_t)  J_{s,t},\ \ \ff{\d}{\d t} \tt J_{s,t}= (A+\nn b_t(X_t))\tt J_{s,t},\ \ J_{s,s}=\tt J_{s,s}=I.
 \end{equation}
 By (H1),   we have
 \beq\label{*D}  \|J_{s,t}^{-1}\|\lor \|J_{s,t}\|\le \e^{\int_s^t K_1(s)}\d s,\ \ t\ge s\ge 0.\end{equation}
 Moreover, since $\e^{At}$ is contractive and $\bar K_1(t):= \|\nn b_t\|_\infty$ is locally bounded in $t$, we have
\beq\label{*D'}   \|\tt J_{s,t}\| \le \e^{\int_s^t \bar K_1(s)}\d s,\ \ t\ge s\ge 0.\end{equation}

\beg{thm}\label{T1.1}  Assume {\rm (H1)} and {\rm (H2)}. If   $\E S(T)^{-\ff 1 2}<\infty$,  then
\beq\label{I} P_T(\nn_{\e^{AT} v} f) = \E\big\{f(X_T)  M_T^v  \big\},\ \   v\in \H,\ f\in C_b^1(\R^d)  \end{equation}   holds for
\beg{align*} M_T^v:=  \ff 1 {S(T)}\bigg(&\int_0^T \big\<\si_t^{-1}\e^{At}  J_{t,T}^{-1} v, \d W_{S(t)}\big\> \\
&+\int_0^T\d S(t) \int_t^T {\rm Tr}\Big\{\si_t^{-1} \e^{At}J_{t,r}^{-1}\big(\nn_{J_{r,T}^{-1}v}B_r\big)(X_r) \tt J_{t,r} \si_t \Big\} \d r\bigg).\end{align*}
\end{thm}

This result extends \cite[Theorem 5.1]{W12} where $S(t)\equiv t$ is considered. When  $\H=\R^d$ is finite-dimensional, we may take $A=0$   so that and Theorem \ref{T1.1} with $\tt J=J$ recovers the main result in \cite{W15}. In this case,
 according to \cite{W12},  the integration by parts formula implies that $P_T$ has a density $p_T(x,y)$ with respect to the Lebesgue measure, which is differentiable in $y$  with
 $$\nn_v\log p_T(x,\cdot)(y)= -\E\Big(M_T^v\Big|X_T(x)=y\Big),\ \ v,x\in\R^d.$$     For nonnegative $f\in \B_b(\R^d)$ and $T>0$, let
$${\rm Ent}_{P_T}(f)= P_T(f\log f)- (P_Tf)\log P_Tf$$ be the relative entropy of $f$ with respect to $P_T$.
Below we present some applications of Theorem \ref{T1.1} for the finite-dimensional case (see also \cite{W15} for the Chinese version).

\beg{cor} \label{C1.2} Assume {\rm (H1)}, {\rm (H2)}, $\H=\R^d, A=0$ and $\E S(T)^{-\ff 1 2}<\infty$.
Let $$\bb(T)=  d T\ll_1(T) \ll_2(T)K_2(T)\e^{3TK_1(T)},\ \ T>0.$$ Then:
\beg{enumerate} \item[$(1)$] For any $T>0$ and $ v\in\R^d$,
\beg{equation*}\beg{split} &\| P_T(\nn_v f)\|_\infty\le |v|\cdot \|f\|_\infty \Big(\ll_2(T)\e^{TK_1(T)} \E S(T)^{-\ff 1 2} + \bb(T)\Big),\ \  f\in C_b^1(\R^d),\\
&  \int_{\R^d} |\nn_v p_T(x,\cdot)|(y) \d y\le |v| \Big(\ll_2(T)\e^{TK_1(T)} \E S(T)^{-\ff 1 2} + \bb(T)\Big),\  \ x\in\R^d.\end{split}\end{equation*}
\item[$(2)$] For any $p>1$, there exists a constant $C(p)\ge 1$ such that for any $T>0$,
\beg{equation*}\beg{split} &| P_T(\nn f)|\le C(p) (P_T|f|^p)^{\ff 1 p} \Big(\ll_2(T)\e^{TK_1(T)} \big(\E S(T)^{-\ff p {2(p-1)}}\big)^{\ff {p-1}p}
  + \bb(T)\Big),\ f\in C_b^1(\R^d),\\
 &  \int_{\R^d} |\nn\log p_T(x,\cdot)|^{\ff p{p-1}}(y)p_T(x,y)\d y \\
  &\qquad  \le C(p) \Big(\ll_2(T)\e^{TK_1(T)} \big(\E S(T)^{-\ff p {2(p-1)}}\big)^{\ff {p-1}p}
  + \bb(T)\Big),\ \ x\in\R^d. \end{split}\end{equation*}
\item[$(3)$] For any $\dd>0,v,x\in\R^d$ and positive $f\in \B_b(\R^d)$,
\beg{equation*}\beg{split} & |P_T(\nn_v f)|\le  \dd{\rm Ent}_{P_T}(f)
+   (P_Tf) \Big(\bb(T)|v|
 + \dd\log\E\exp\Big[\ff{\ll_2(T)^2|v|^2\e^{2TK_1(T)}}{2\dd^2 S(T)}\Big]\Big),\\
 &  \int_{\R^d} \exp\Big[\ff{|\nn_v\log p_T(x,\cdot)|(y)}\dd\Big]p_T(x,y)\d y
 \le \E\exp\Big[\ff{\bb(T)|v|}\dd +\ff{\ll_2(T)^2|v|^2\e^{2TK_1(T)}}{2\dd^2S(T)} \Big].\end{split}\end{equation*} \end{enumerate}
\end{cor}

\beg{cor} \label{C1.3} In the situation of Corollary  $\ref{C1.2}$.    Let $p>1, T>0.$ If
$$\GG_{T,p}(r):= \E\exp\bigg[\ff{p^2\ll_2(T)^2\e^{2TK_1(T)}r^2}{2(p-1)^2S(T)}\bigg]<\infty,\ \   r\ge 0,$$ then the shift Harnack inequality
\beq\label{SH}(P_Tf)^p(x)\le \exp\Big[ \ff{p(\log p)\bb(T)|v|}{p-1}
+\ff{p-1}p\log \GG_{T,p}(|v|)\Big]P_T(f^p(v+\cdot))(x)\end{equation} holds for all $v,x\in\R^d$ and positive $f\in \B_b(\R^d).$ Consequently,
$$\sup_{x\in\R^d} \int_{\R^d} p_T(x,y)^{\ff p{p-1}}\d y \le \bigg(\int_{\R^d}\exp\Big[ -\ff{p(\log p)\bb(T)|v|}{p-1}
-\ff{p-1}p\log \GG_{T,p}(|v|)\Big]\d v\bigg)^{\ff{-1}{p-1}}.$$ \end{cor}

\

To illustrate the above results,   we consider below the SDE driven by $\aa$-stable like noises.

\

\beg{cor}\label{C1.4} In the situation of Corollary $\ref{C1.2}$.  Let $B(r)\ge cr^{\ff\aa 2}$ for $r\ge r_0$, where  $\aa\in (0,2)$ and  $c,r_0>0$ are constants.
\beg{enumerate} \item[$(1)$] For any $p>1$ there exists a constant $C(p)>0$ such that
\beg{equation*}\beg{split} &|P_T(\nn f)|\le \ff{C(p)(P_T|f|^p)^{\ff 1 p}}{1\land T^{\ff 1 \aa}},\ \ T>0, f\in C_b^1(\R^d),\\
 &\sup_{x\in\R^d} \int_{\R^d} |\nn\log p_T(x,\cdot)|^{\ff p{p-1}}(y)p_T(x,y)\d y\le \ff{C(p)}{1\land T^{\ff 1 \aa}},\ \ T>0.\end{split}\end{equation*} \item[$(2)$] Let $\aa\in (1,2).$ Then there exists a constant $C>0$ such that for any $p>1,\dd>0, v\in\R^d$ and $f\in C^1(\R^d)$,
 \beg{equation*}\beg{split}&|P_T(\nn_v f)|\le \dd{\rm Ent}_{P_T}(f) + (P_Tf) \Big(\bb(T)|v|+\ff{C|v|^2}{\dd^2(1\land T)^{\ff 2\aa}}+ \ff{C|v|^{\ff\aa{\aa-1}} }{\{\dd^\aa (1\land T)\}^{\ff 1 {\aa-1}}}\Big),\\
 &\sup_{x\in\R^d}\int_{\R^d} \exp\Big[\ff{|\nn_v\log p_T(x,\cdot)(y)|}{\dd}\Big]p_T(x,y)\d y\\
 &\qquad \le \exp\Big[\bb(T)|v|+ \ff{C|v|^2}{\dd^2(1\land T)^{\ff 2\aa}}+ \ff{C|v|^{\ff\aa{\aa-1}} }{\{\dd^\aa (1\land T)\}^{\ff 1 {\aa-1}}}\Big].\end{split}\end{equation*}
 \item[$(3)$] Let $\aa\in (1,2).$ Then there exists a constant $C>0$ such that for any $p>1,T>0, v\in\R^d$ and positive $f\in\B_b(\R^d), $
 \beg{equation*}\beg{split} &(P_Tf)^p\le \exp\Big[\ff{C(p\log p)|v|}{p-1}+\ff{Cp|v|^2}{(p-1)(1\land T)^{\ff 2\aa}}+ \ff{Cp^{\ff 1{\aa-1}}|v|^{\ff\aa{\aa-1}}}{[(p-1)(1\land T)]^{\ff 1{\aa-1}}}\Big]P_T( f^p(v+\cdot)),\\
 &\sup_{x\in\R^d} \int_{\R^d} p_T(x,y)^{\ff p{p-1}}\d y\le \ff{1}{(1\land T)^{\ff d{\aa(p-1)}}}\exp\Big[\ff{Cp\log p}{(p-1)^2}+\ff{Cp^{\ff 1 {\aa-1}}}{(p-1)^{\ff\aa{\aa-1}}}\Big].\end{split}\end{equation*} \end{enumerate}\end{cor}

 The remainder of the paper is organized as follows. In Section 2, we fix a path $\ell$ of $S$ with finite jumps, and establish  the integration by parts formula for the corresponding equation, i.e. the equation \eqref{1.1} with $\ell$ in place of $S$. In Section 3, we use this integration by parts formula to prove   the above results by using finite jump approximations.

 \section{Integration by parts formula for the equation  with finite jump}

 In this section, we let $\ell$ be a c\'adl\'ag and increasing function on $[0,\infty)$ with $\ell(0)=0$ such that the set $\{t\in [0,T]: \DD \ell(t):= \ell(t)-\ell(t-)>0\}$ is finite. We call $\ell $ a path of $S$ with finite many jumps on $[0,T]$. Let $X_t^\ell$ solve the equation
\beq\label{2.10}   X_t^\ell =\e^{At} X_0^\ell + \int_0^t \e^{A(t-s)}b_s(X_s^\ell)\d s+ \int_0^t\e^{A(t-s)}\si_s\d W_{\ell(s)}+V_t,\ \ t\ge 0,\end{equation}   and let $P_t^\ell$ be the associated Markov operator; i.e.
 $$P_t^\ell f(x):= \E f(X_t^\ell(x)),$$ where $X_t^\ell(x)$ solves \eqref{2.10} for $X_0^\ell=x.$ Moreover, let $(J_{s,t}^\ell, \tt J_{s,t}^\ell)_{t\ge s}$ be defined  in \eqref{*0} for $X^\ell$ in place of $X$.
 The main result in this section is the following.

 \beg{thm}\label{T2.1} Assume {\rm (H1)} and {\rm (H2)}. $\ell$ be a path of $S$ with finite many jumps on $[0,T]$  and $\ell(T)>0$.   Then
\beq\label{I} P_T^\ell(\nn_{\e^{AT} v} f) = \E\big\{f(X_T)  M_T^{\ell,v} \big\},\ \   v\in \H,\ f\in C_b^1(\R^d)  \end{equation}   holds for
\beg{align*} M_T^{\ell,v}:=  \ff 1 {\ell(T)}\bigg(&\int_0^T \big\<\si_t^{-1}\e^{At}  (J_{t,T}^\ell)^{-1} v, \d W_{\ell(t)}\big\> \\
&+\int_0^T\d \ell(t) \int_t^T {\rm Tr}\Big\{\si_t^{-1} \e^{At}(J_{t,r}^\ell)^{-1}\big(\nn_{(J_{r,T}^\ell)^{-1}v}B_r\big)(X_r^\ell) \tt J_{t,r}^\ell \si_t \Big\} \d r\bigg).\end{align*}   \end{thm}

\beg{proof} We shall use the integration by parts formula in the Malliavin calculus, see, for instance  \cite{M,N}. For the Brwonian motion $(W_t)_{t\in [0,\ell(T)]},$ let $(D,\D(D))$ be the    Malliavin gradient, and let $(\ddd,\D(\ddd))$ be its adjoint operator (i.e. the Malliavin divergence).    Let $J_t^\ell= J_{0,t}^\ell$ and $\tt J_t^\ell= \tt J_{0,t}^\ell$. It is easy to see that
\beq\label{JJ}  J_T^\ell =  J_{t,T}^\ell   J_t^\ell,\  \tt J_T^\ell = \tt J_{t,T}^\ell \tt J_t^\ell,\    \tt J_t^\ell= \e^{At} J_t^\ell,\ \ T\ge t\ge 0.\end{equation} Take
$$h(t)= \sum_{i=1}^t \big(t\land \ell(t_i)-\ell(t_{i-1})\big)^+\si_{t_i}^{-1}\tt J^\ell_{t_i}  (J_T^\ell )^{-1} v,\ \ t\in [0,\ell(T)].$$ From  (H1) we see that $J_t^\ell$ and $ (J_t^\ell)^{-1}$ are Malliavin differentiable for every $t\in [0,\ell(T)],$  such that $h\in \D(\ddd)$.
Since $(V_t)_{t\ge 0}$ is independent of $(W_t)_{t\ge 0}$, we have $D_hV_t=0$, so that (\ref{2.10}) yields
\beq\label{D*} \d D_h X^\ell_t =\big\{A+  (\nn b_t)(X^\ell_t)\big\} D_h X^\ell_t \,\d t +\si_t \d h_{\ell(t)},\ \ D_h X^\ell_0=0.\end{equation} Then
by Duhamel's formula and \eqref{JJ},   $$D_h X^\ell_T=  \int_0^T\tt J^\ell_{t,T} \si_t \d h_{\ell(t)} =  \sum_{i=1}^n  \tt J^\ell_{t_i,T}  \si_{t_i}\si_{t_i}^{-1}\tt J^\ell_{t_i} \DD \ell(t_i) (J_T^\ell )^{-1} v =\ell(T)\e^{AT} v.$$
Therefore, \beq\beg{split} \label{2.1}  &\E(\nn_{\e^{AT} v} f)(X^\ell_T) =\ff 1 {\ell(T)} \E (\nn_{D_hX^\ell_T} f)(X^\ell_T) \\
&= \ff 1 {\ell(T)} \E\big\{D_h f(X^\ell_T)\big\} = \ff 1 {\ell(T)} \E\big\{f(X^\ell_T) \ddd(h)\big\}.\end{split}\end{equation}

To calculate $\ddd(h)$, let
$$h_{ik}(t)= (t\land \ell(t_i)-\ell(t_{i-1}))^+ e_k,\ \ F_{ik} =\big\<\si_{t_i}^{-1} \tt J^\ell_{t_i} (J^\ell_T)^{-1} v, e_k\big\>$$ for $1\le i\le n, 1\le k\le d, t\in [0,\ell(T)],$ where $\{e_k\}_{k=1}^d$ is the canonical orthonormal basis on $\R^d$. Then
$$h(t)= \sum_{k=1}^d\sum_{i=1}^n F_{ik} h_{ik}(t),\ \ t\in [0,\ell(T)].$$
Noting that $h_{ik}$ is deterministic with $\int_0^{\ell(T)} |h'_{ik}(t)|^2\d t<\infty$, we have  $$\ddd(h_{ik})= \int_0^{\ell(T)} \<h_{ik}'(t),\d W_t\>= \<e_k, W_{\ell(t_i)}-W_{\ell(t_{i-1})}\>.$$   Thus, using the formula $\ddd(F_{ik}h_{ik})= F_{ik}\ddd(h_{ik})-D_{h_{ik}}F_{ik},$ we obtain
\beq\beg{split}\label{2.2} \ddd(h) &= \sum_{k=1}^d\sum_{i=1}^n \big\{F_{ik} \ddd(h_{ik}) -D_{h_{ik}}F_{ik}\big\}\\
&=\sum_{k=1}^d\sum_{i=1}^n \big\{F_{ik} \<e_k, W_{\ell(t_i)}-W_{\ell(t_{i-1})}\> -\<\si_{t_i}^{-1} D_{h_{ik}}(\tt J^\ell_{t_i}(J^\ell_T)^{-1})v,e_k\>\big\}\\
&=  \int_0^T\big\<\si_t^{-1} \tt J^\ell_t(J^\ell_T)^{-1} v,   \d W_{\ell(t)}  \big\> - \sum_{k=1}^d\sum_{i=1}^n \big\<\si_{t_i}^{-1} D_{h_{ik}}(\tt J^\ell_{t_i}(J^\ell_T)^{-1})v, e_k\big\>.\end{split}\end{equation}  Since   $\d h_{ik} (t)$ is supported on $(\ell(t_{i-1}), \ell(t_i))$ but $\tt J^\ell_{t_i}$ is determined by
$ \{W_{t}: t\le \ell(t_{i-1}) \}$, we have $D_{h_{ik}}\tt J_{t_i}^\ell =0$ so that
\beq\label{2.3} D_{h_{ik}} (\tt J^\ell_{t_i}(J^\ell_T)^{-1}) = \tt J^\ell_{t_i} D_{h_{ik}} (J^\ell_T)^{-1}= -\tt J^\ell_{t_i} (J^\ell_T)^{-1} (D_{h_{ik}} J^\ell_T) (J^\ell_T)^{-1}.\end{equation} Noting that  (\ref{*0}) for $J_t^\ell:= J_{0,t}^\ell$  yields
$$\d D_{h_{ik}} J^\ell_t= (\nn_{D_{h_{ik}}X^\ell_t} B_t)(X^\ell_t)J^\ell_t \d t +  B_t (X^\ell_t) D_{h_{ik}}J^\ell_t \d t,\ \ D_{h_{ik}}J^\ell_0=0,$$ by Duhamel's formula  we obtain
\beq\label{2.4}D_{h_{ik}}J^\ell_T =   \int_0^T  J^\ell_{t,T}   (\nn_{D_{h_{ik}} X^\ell_t} B_t)(X^\ell_t) J^\ell_t\d t.\end{equation} Moreover, it follows from (\ref{D*}) that
$$D_{h_{ik}} X^\ell_t=   \int_0^t \tt J^\ell_{s,t}  \si_s \d (h_{ik})_{\ell(s)} = 1_{\{t_i\le t\}} (\ell(t_i)-\ell(t_{i-1})) \tt J^\ell_{t_i,t} \si_{t_i} e_k.$$
Combining this with (\ref{2.4}),  we arrive at
$$D_{h_{ik}} J^\ell_T= (\DD\ell(t_i))  \int_{t_i}^T  J^\ell_{s,T}  \big(\nn_{\tt J^\ell_{t_i,s}  \si_{t_i} e_k}B_s\big)(X^\ell_s) J^\ell_s \d s.$$ Substituting this into (\ref{2.3}) we obtain
\beg{equation*}\beg{split} &\sum_{k=1}^d\sum_{i=1}^n \big\<\si_{t_i}^{-1} D_{h_{ik}} (\tt J^\ell_{t_i}(J^\ell_T)^{-1})v,e_k\big\>\\
&= -\sum_{k=1}^d\sum_{i=1}^n (\DD\ell(t_i))  \int_{t_i}^T \big\< \si_{t_i}^{-1} \tt J^\ell_{t_i}(J^\ell_r)^{-1} \big(\nn_{\tt J^\ell_{t_i,r}  \si_{t_i} e_k}B_r\big)(X^\ell_s)(J^\ell_{r,T})^{-1}v,e_k\big\>\d r \\
&=- \sum_{k=1}^d\int_0^T \d\ell(t) \int_t^T \Big\<\si_t^{-1}\tt J^\ell_t (J^\ell_r)^{-1} \big(\nn_{(J^\ell_{r,T})^{-1}v} B_r\big)(X_r^\ell) \tt J_{t,r}^\ell \si_t e_k, e_k\Big\>\d r\\
&=-\int_0^T\d\ell(t)\int_t^T{\rm Tr}\Big\{  \si_t^{-1}\tt J^\ell_t (J^\ell_r)^{-1} \big(\nn_{(J^\ell_{r,T})^{-1}v}B_r\big)(X^\ell_r)\tt
J^\ell_{t,r}  \si_t  \Big\}\d s.\end{split}\end{equation*}
Therefore, we derive from \eqref{JJ} and (\ref{2.2})  that $\ddd(h)= M_T^{\ell,v}$ and hence,      the proof is finished by (\ref{2.1}).\end{proof}

\section{Proofs}

\beg{proof}[Proof of Theorem \ref{T1.1}] According to \cite[Theorem 2.4(1)]{W12}, the second assertion follows from the first. So, it suffices to prove the desired integration by parts formula.  For any path $\ell$ of $S$ with $\ell(T)>0$, for any $\vv>0$, let
$$\ell_\vv(t)= \sum_{s\le t} \DD\ell(s) 1_{\{\DD\ell(s) \ge\vv\}},\ \ t\ge 0.$$ Then $\ell_\vv$ has finite many jumps on $[0,T].$ Moreover, $\d\ell_\vv(t)\to
\d\ell(t)$ on $[0,T]$ strongly as $\vv\to 0.$ Note that by \eqref{*D'}, (H1) and (H2),
$$\big\|\si_t^{-1}(J_{t,T}^{\ell_\vv})^{-1}\big\|+
 \int_t^T \big\|\si_t^{-1} (J_{t,r}^{\ell_\vv})^{-1}   \big(\nn_{(J_{t,r}^{\ell_\vv})^{-1}\si_t e_k}B_r\big)(X_r^{\ell_\vv})
 \tt J_{t,r}^{\ell_\vv} \big\| \d r$$ is  bounded in $(t,\vv)\in [0,T]\times [0,1],$ and by \cite[Lemma 3.1]{WXZ}
$$\lim_{\vv\to 0} \E \sup_{t\in [0,T]} |X_t^{\ell_\vv} -X^\ell_t|^2 =0,$$ which together with (H1) implies
$$\lim_{\vv\to 0} \E \sup_{t\in [0,T]} \big(\|J_t^{\ell_\vv} -J_t^\ell\|^2+ \|(J_t^{\ell_\vv})^{-1}- (J_t^\ell)^{-1}\|^2\big)=0.$$
Due the contraction of $\e^{At}$ and the second formula in \eqref{JJ}, the same holds for $\tt J_t $ in place of $J_t$. Combining these with  \eqref{*D'}, (H1) and (H2),  we conclude that
\beg{equation*}\beg{split} &\lim_{\vv\to 0} P_T^{\ell_\vv} (\nn_vf) = \lim_{\vv\to 0} \E(\nn_vf)(X_T^{\ell_\vv}) = P_T^\ell (\nn_vf),\\
&\lim_{\vv\to 0} \E\big\{f(X_T^{\ell_\vv})M_T^{\ell_\vv,v}\big\}= \E\big\{f(X_T^{\ell})M_T^{\ell,v}\big\},\ \ f\in C_b^1(\R^d).\end{split}\end{equation*}
Therefore, first applying Theorem \ref{T2.1} to $\ell_\vv$ in place of $\ell$ then letting $\vv\downarrow 0$, we obtain
$$P_T^\ell(\nn_v f)= \E\big\{f(X_T^\ell)M_T^{\ell,v}\big\}$$ for all sample path $\ell$ of $S$ with $\ell(T)>0$. Since
$\E S(T)^{-\ff 1 2}<\infty$ implies $S(T)>0$, and noting that $X_T= X_T^S, M_T^v=M_T^{S,v}$, we obtain
\beq\label{D3} P_T^S(\nn_vf) =   \E^S\big\{f(X_T)M_T^{v}\big\},\end{equation}where   $\E^S$ is the conditional expectation given $S$.
Moreover,  it follows from   \eqref{*D'}, (H1), (H2),   and $\E S(T)^{-\ff 1 2}<\infty$ that
\beq\label{D*2} \beg{split} &\E\big| M_T^v  \big|= \E\big[ \E^S |M_T^v|\big] \\
&\le  \E \bigg[\ff 1 {S(T)}  \bigg(\E^S \int_0^T |\si_t^{-1}J_{t,T}^{-1}v|^2\d S(t)\bigg)^{1/2}\\
&\quad + \sum_{k=1}^d \ff 1 {S(T)}\int_0^T \d S(t) \int_t^T \|\si_t^{-1} J_{t,r}^{-1}\|\cdot\big| \big(\nn_{J_{t,r}^{-1}v} B_r)(X_r) \tt J_{t,r} \si_t e_k\big|\d r\bigg]\\
&\le |v| \Big(\ll_2(T)\e^{TK_1(T)}\E S(T)^{-\ff 1 2} + d T \ll_1(T)\ll_2(T) K_2(T) \e^{3TK_1(T)}\Big)<\infty.\end{split}\end{equation}   Then $  M_T^v\in L^1( \P)$ so that  (\ref{D3}) yields
$$P_T(\nn_vf)= \E P_T^S(\nn_vf) = \E\big[ f(X_T) M_T^v\big].$$ This completes the proof.
\end{proof}

\beg{proof}[Proof of Corollary \ref{C1.2}] We have $\tt J_t= J_t$. Assertion (1) follows immediately from (\ref{D*2}),   Theorem \ref{T1.1} and \cite[Theorem 2.4(1)]{W12} with $H(r)=r$.

Next,    by  \eqref{*D}, (H1), (H2)  and the Burkholder inequality \cite[Theorem 2.3]{Z} (see also \cite[Lemma 2.1]{WXZ}),  for any $p>1$ there exists a constant $C(p)\ge 1$ such that
\beg{equation*}\beg{split}    \big(\E  |M_T^v|^{\ff p{p-1}} \big)^{\ff {p-1}p}
&\le   \bb(T)|v|+C(p) \Big(\E \ff {(\int_0^T|\si_t^{-1} J_{t,T}^{-1}v|^2\d S(t))^{\ff p{2(p-1)}}}  {S(T)^{\ff p{p-1}}} \Big)^{\ff{p-1}p}\\
&\le \bb(T)|v|+ C(p)|v|   \ll_2(T)\e^{TK_1(T)}\big(\E S(T)^{-\ff p {2(p-1)}}\big)^{\ff {p-1}p}. \end{split}\end{equation*}    Then assertion (2) follows from \cite[Theorem 2.4(1)]{W12} with $H(r)= r^{\ff{p}{p-1}}$ and the fact that
$$|P_T(\nn_vf)|= \Big|\E\big\{f(X_T) M_T^v \big\}\Big|\le (P_T|f|^p)^{\ff 1 p} \big(\E  |M_T^v|^{\ff p{p-1}} \big)^{\ff {p-1}p},\ \ v\in\R^d.$$

Finally, by Theorem \ref{T1.1} and the Young inequality (see \cite[Lemma 2.4]{ATW09}), if $f\in C_b^1(\R^d)$ is nonnegative, then
\beq\label{2*} |P_T(\nn_vf)|= \Big|\E\big\{ f(X^\ell_T)M_T^v \big\}\Big| \le\dd {\rm Ent_{P_T}}(f) +\dd (P_T f) \log \E\exp\Big[ \ff{M_T^v}{\dd }\Big],\ \dd>0.\end{equation} Obviously, by  \eqref{*D}, (H1) and (H2),
\beg{equation*}\beg{split} & M_T^v \le  \bb(T)|v| +\ff 1 {S(T)} \int_0^T \<\si_t^{-1}J_{t,T}^{-1}v, \d W_{S(t)}\>,\\
& \E^S \exp\bigg[\ff 1 {\dd S(T)} \int_0^T \<\si_t^{-1}J_{t,T}^{-1}v, \d W_{S(t)}\>\bigg] \le \exp\bigg[\ff{\ll_2(T)^2|v|^2\e^{2TK_1(T)}}{2\dd^2S(T)}\bigg],\ \ \dd>0.\end{split}\end{equation*} Then
$$\log \E\exp\Big[\ff{M_T^v}{\dd }\Big]
 \le \ff{\bb(T)|v|}\dd + \log \E \exp\bigg[\ff{\ll_2(T)^2|v|^2\e^{2TK_1(T)}}{2\dd^2 S(T)} \bigg],\ \ \dd>0.$$ By combining this with (\ref{2*}) and \cite[Theorem 2.4(1)]{W12} for $H(r)= \e^{r/\dd}$, we prove (3).
\end{proof}

\beg{proof}[Proof of Corollary \ref{C1.3}] By \cite[Theorem 2.5(2)]{W12}, the second assertion follows from the first. So, we only need to prove the required shift Harnack inequality \eqref{SH} for $v\ne 0$.  By Corollary \ref{C1.2}(3), we have
 $$|P_T(\nn_v f)|\le \dd {\rm Ent}_{P_T}(f) +   (P_T f)\Big(\bb(T)|v| + \dd \log \E\exp\Big[\ff{\ll_2(T)|v|^2}{2\dd^2 S(T)}\e^{\int_0^TK_1(t)\d t}\Big]\Big),\ \ \dd>0.$$
So, letting
$$\bb_v(\dd)= \bb(T)|v| + \dd \log \E\exp\Big[\ff{\ll_2(T)^2|v|^2\e^{2TK_1(T)}}{2\dd^2 S(T)}\Big],\ \ \dd>0,$$
we obtain from \cite[Proposition 2.3]{W12} that
\beq\label{*W} (P_Tf)^p\le (P_Tf^p(v+\cdot)) \exp\bigg[\int_0^1\ff{p}{1+(p-1)s} \bb_v\Big(\ff{p-1}{1+(p-1)s}\Big)\d s\bigg].\end{equation}
By the Jensen inequality,  for $\dd= \ff{p-1}{1+(p-1)s}$ we have
\beg{equation*}\beg{split} \E \exp\Big[\ff{\ll_2(T)^2|v|^2\e^{2TK_1(T)}}{2\dd^2 S(T)}\Big]&\le \bigg(\E\exp\Big[\ff{p^2\ll_2(T)^2|v|^2\e^{2TK_1(T)}}{2(p-1)^2S(T)}\Big]\bigg)^{\ff{(1+(p-1)s)^2}{p^2}}\\
&=\GG_{T,p}(|v|)^{\ff{(1+(p-1)s)^2}{p^2}}.\end{split}\end{equation*} Thus,
\beg{equation*}\beg{split} &\int_0^1\ff{p}{1+(p-1)s}\bb_v\Big(\ff{p-1}{1+(p-1)s}\Big)\d s\\
 &\le \bb(T)|v|\int_0^1\ff p{1+(p-1)s}\d s +\ff{p-1}p\log \GG_{T,p}(|v|)\d s\\
&= \ff{p\log p}{p-1} \bb(T)|v| + \ff{p-1}p\log\GG_{T,p}(|v|).\end{split}\end{equation*}
Then the proof is finished by combining this with (\ref{*W}).
\end{proof}

\beg{proof}[Proof of Corollary \ref{C1.4}] Since assertions in Corollaries \ref{C1.2} and \ref{C1.3} are unform in $V$, we may apply them for any deterministic path of $V$ in place of the process $V$, so that these two Corollaries remain true for $P_T^V$ in place of $P_T$, where
$$P_T^Vf(x) = \E^V (f(X_T(x)):= \E\big(f(X_T(x))|V\big).$$

Next, we observe that by the Markov property it suffices to prove the assertions for  $P_T^V$ in place of $P_T$ with $T\in (0,1]$.  In fact, for $T>1$ let
$$P_{1,T}^V f(x)= \E^V f(X_{1,T}(x)),\ \ f\in\B_b(\R^d), x\in\R^d,$$ where $(X_{1,t}(x))_{t\ge 1}$ solves the equation
$$X_{1,t}(x)= x + \int_1^tb_s(X_{1,s}(s))\d s +\int_1^t \si_s\d W_{S(s)}+V_t-V_1,\ \ t\ge 1.$$ Then by the Markov property of $X_t$ under $\E^V$, we obtain,
$$P_T^Vf= P_{1,T}^V(P_1^V f),\ \ f\in\B_b(\R^d).$$ Combining this with the assertions for $T=1$ and using the Jensen inequality, we prove the assertions for $T>1$. For instance, if for $p>1$ one has
$$|P_1^V(\nn f)|\le C(p) (P_1^V|f|^p)^{\ff 1 p},$$ then for any $T>1,$  \beg{equation*}\beg{split} &|P_T(\nn f)|= |\E P_{1,T}^VP_1^V(\nn f)|\le \E P_{1,T}^V|P_1^V(\nn f)|\\
&\le C(p) \E P_{1,T}^V (P_1^V|f|^p)^{\ff 1 p}\le C(p)(P_T|f|^p)^{\ff 1 p}= \ff{C(p)(P_T|f|^p)^{\ff 1 p}}{(1\land T)^{\ff  1 \aa}}.\end{split}\end{equation*} Below we prove assertions (1)-(3) for $T\in (0,1]$ respectively.

(1) Since $\bb(T)+ \ll_2(T)\e^{TK_1(T)}$ is bounded for $T\in (0,1]$, and by \cite[(ii) in the proof of Theorem 1.1]{WXZ}
$$\Big(\E S(T)^{-\ff{p}{2(p-1)}}\Big)^{\ff{p-1}p}\le \ff C{T^{\ff 1 \aa}},\ \ T\in (0,1]$$ holds for some constant $C>0$, the desired assertion follows from Corollary \ref{C1.2}(2).

(2) Let $\aa\in (1,2),$  and let $S_\aa$ be the subordinator induced by  the Bernstein function $r\mapsto r^{\ff\aa 2}$. Then as shown in \cite[Proof of Corollary 1.2]{WW} that
$$\E \ff 1 {S(T)^k}\le c_0\E\ff 1 {S_\aa(T)^k},\ \ k\ge 1, T\in (0,1]$$ holds for some constant $c_0\ge 1.$ Combining this with the third display from below  in the proof of \cite[Theorem 1.1]{GRW} for
$\kk=1$, i.e. (note the $\aa$ therein is $\aa/2$ here)
$$\E \e^{\ll/\tt S(t)} \le 1+\bigg(\exp\Big[\ff{c_1\ll^{\ff{\aa}{2(\aa-1)}}}{t^{\ff 1 {\aa-1}}}\Big]-1\bigg)^{\ff{2(\aa-1)}\aa}\le
\exp\Big[ \ff{c_2\ll}{t^{\ff 2 \aa}}+\ff{c_2\ll^{\ff\aa{2(\aa-1)}}}{t^{\ff 1{\aa-1}}}\Big],\ \ \ll,t\ge 0$$ for some constants $c_1,c_2>0,$ we obtain
\beq\label{*F}\beg{split}  \E\e^{\ll/S(T)} &\le 1+c_0\big(\E\e^{\ll/S_\aa(T)}-1\big)\le  \E\e^{c_0\ll/S_\aa(T)} \\
&\le \exp\Big[ \ff{c_3\ll}{t^{\ff 2 \aa}}+\ff{c_3\ll^{\ff\aa{2(\aa-1)}}}{t^{\ff 1{\aa-1}}}\Big],\ \ T\in (0,1], \ll\ge 0 \end{split}\end{equation} for some constant $c_3>0$.
By Corollary \ref{C1.2}(3) and (\ref{*F}), we prove the desired assertion.

(3) By (\ref{*F}), there exists a constant $c_4>0$ such that
$$\GG_{T,p}(r) \le \exp\Big[\ff{c_4 p^2r^2}{(p-1)^2T^{\ff 2\aa}}+\ff{c_4(pr)^{\ff\aa{\aa-1}}}{(p-1)^{\ff\aa{\aa-1}}T^{\ff 1{\aa-1}}}\Big],\ \ \ r\ge 0, T\in (0,1].$$ Then  there exists a constant $c_5>0$ such that
\beq\label{*F2} \beg{split}&\ff{p(\log p)\bb(T)|v|}{p-1}+ \ff{p-1}p \log \GG_{T,p}(|v|)\\
&\le \ff{c_5(p\log p)|v|}{p-1}
+ \ff{c_5 p|v|^2}{(p-1)T^{\ff 2\aa}}+\ff{c_5p^{\ff 1 {\aa-1}}|v|^{\ff\aa{\aa-1}}}{(p-1)^{\ff1{\aa-1}}T^{\ff 1{\aa-1}}},\ \ T\in (0,1], v\in\R^d.\end{split}\end{equation} By Corollary \ref{C1.3}, this implies the first inequality in (3) for some constant $C>0$.
Finally, the second inequality in (3) follows since (\ref{*F2}) and Corollary \ref{C1.3} imply
\beg{equation*}\beg{split} \sup_{x\in\R^d} \int_{\R^d} p_T(x,y)^{\ff p{p-1}}\d y &\le \bigg(\int_{\R^d} \exp\Big[-\ff{C(p\log p)|v|}{p-1}
-\ff{Cp|v|^2}{(p-1)T^{\ff 2\aa}}-\ff{Cp^{\ff 1{\aa-1}}|v|^{\ff\aa{\aa-1}}}{[(p-1)T]^{\ff 1{\aa-1}}}\Big]\d v\bigg)^{\ff {-1}{p-1}}\\
&\le \bigg(\int_{\{|v|\le T^{^\ff 1 \aa}\}} \d v\bigg)^{\ff{-1}{p-1}}\exp\Big[\ff{Cp(1+\log p)}{(p-1)^2}+\ff{Cp^{\ff 1 {\aa-1}}}{(p-1)^{\ff\aa{\aa-1}}}\Big]\\
&\le \ff 1 {T^{\ff d{\aa(p-1)}}}\exp\Big[\ff{C'p\log p}{(p-1)^2}+\ff{C'p^{\ff 1 {\aa-1}}}{(p-1)^{\ff\aa{\aa-1}}}\Big]\end{split}\end{equation*} for some constant $C'\ge C$.
\end{proof}

\paragraph{Acknowledgement.} The author would like to thank professor Xicheng Zhang and the referee for helpful comments.

\end{document}